\documentclass{article}
\usepackage{amsmath}
\usepackage{amsfonts}
\usepackage{amsthm}
\usepackage{amssymb}
\usepackage{color}

\newtheorem{theorem}{Theorem}[section]

\theoremstyle{definition}
\newtheorem{definition}[theorem]{Definition}
\theoremstyle{remark}
\newtheorem{remark}[theorem]{Remark}

\def\f2{\mathbb{F}_2}

\newcommand{\RNP}{{\rm RNP}}

\begin{document}

\title{Metric spaces nonembeddable into Banach spaces with the Radon-Nikod\'ym property and thick families of geodesics}

\author{Mikhail Ostrovskii\footnote{Supported in part by NSF
DMS-1201269.}{~}\footnote{The author would like to thank Sean Li
for useful information and the referee for helpful criticism.}}

\date{\today}
\maketitle

\begin{large}

\noindent{\bf Abstract.} We show that a geodesic metric space
which does not admit bilipschitz embeddings into Banach spaces
with the Radon-Nikod\'ym property does not necessarily contain a
bilipschitz image of a thick family of geodesics. This is done by
showing that any thick family of geodesics is not Markov convex,
and comparing this result with results of Cheeger-Kleiner,
Lee-Naor, and Li. The result contrasts with the earlier result of
the author that any Banach space without the Radon-Nikod\'ym
property contains a bilipschitz image of a thick family of
geodesics.
\medskip

\medskip

\medskip

\section{Introduction}

The Radon-Nikod\'ym property (\RNP) is one of the most important
isomorphic invariants of Banach spaces. We refer to \cite{BL00,
Bou79, Bou83, DU77, Dul85, Pis11} for systematic presentations of
results on the \RNP.\medskip

In the recent work on metric embeddings a substantial role is
played by existence and non-existence of bilipschitz embeddings of
metric spaces into Banach spaces with the \RNP, see
\cite{CK06,CK09,LN06}. At the seminar ``Nonlinear geometry of
Banach spaces'' (Texas A \&\ M University, August 2009) Bill
Johnson suggested the problem of metric characterization of
reflexivity and the \RNP~ \cite[Problem 1.1]{Tex09}; see also
\cite[p.~307]{Ost13}. In \cite{Ost14} the RNP was characterized in
terms of thick families of geodesics defined in the following way:

\begin{definition}\label{D:ThickFam} Let $u$ and $v$ be two elements in a metric space $(M,d_M)$. A {\it $uv$-geodesic} is a
distance-preserving map $g:[0,d_M(u,v)]\to M$ such that $g(0)=u$
and $g(d_M(u,v))=v$ (where $[0,d_M(u,v)]$ is an interval of the
real line with the distance inherited from $\mathbb{R}$). A family
$T$ of $uv$-geodesics is called {\it thick} if there is $\alpha>0$
such that for every $g\in T$ and for any finite collection of
points $r_1,\dots,r_n\in [0,d_M(u,v)]$, we are going to call them
{\it control points}, there is another $uv$-geodesic $\widetilde
g\in T$ and a sequence $0< s_1< q_1< s_2< q_2<\dots< s_m<q_m<
s_{m+1}<d_M(u,v)$ satisfying the conditions:

\begin{itemize}

\item The set $\{0,q_1,\dots,q_m,d_M(u,v)\}$ contains
$r_1,\dots,r_n$.

\item $g(q_i)=\widetilde g(q_i)$.

\item
\begin{equation}\label{E:alphaSep}\sum_{i=1}^{m+1}d_M(g(s_i),\widetilde
g(s_i))\ge\alpha.\end{equation}

\end{itemize}

\end{definition}

The following result gives a metric characterization of the RNP.

\begin{theorem}[\cite{Ost14}]\label{RNPThick}
A Banach space $X$ does not have the \RNP\ if and only if  there
exists a metric space $M_X$ containing a thick family $T_X$ of
geodesics which admits a bilipschitz embedding into $X$.
\end{theorem}

Studying metric characterizations of the RNP, it would be much
more useful and interesting to get a characterization of all
metric spaces which do not admit bilipschitz embeddings into
Banach spaces with the RNP. In view of Theorem \ref{RNPThick} it
is natural to ask whether the presence of bilipschitz images of
thick families of geodesics characterizes metric spaces which do
not admit bilipschitz embeddings into spaces with the RNP? It is
clear that the answer to this question in full generality is
negative: we may just consider a dense subset of a Banach space
without the RNP which does not contain any continuous curves (e.
g. subset of all vectors with rational coordinates in $c_0$). So
we need to restrict our attention to spaces containing
sufficiently large collections of continuous curves. Our main
result is a negative answer even in the case of geodesic metric
spaces (we use the terminology of \cite{BH99} on metric spaces and
of \cite{Ost13} on metric embeddings):

\begin{theorem}\label{T:main} There exist geodesic metric spaces which
\begin{itemize}

\item Do not contain bilipschitz images of thick families of
geodesics.

\item Do not admit bilipschitz embeddings into Banach spaces with
the Radon-Niko\-d\'ym property.
\end{itemize}
\end{theorem}

More precisely, we prove that the Heisenberg group with its
subriemannian (Carnot-Caratheodory) metric (see \cite{CDPT07,
Gro96, Li13, Mon02}) does not admit a bilipschitz embedding of a
thick family of geodesics. This result proves Theorem \ref{T:main}
since it is known that the Heisenberg group does not admit a
bilipschitz embedding into a Banach space with the RNP, see
\cite{CK06} and \cite{LN06}, where the observation made in
\cite{Sem96} on the consequences of the differentiability result
of \cite{Pan89} was generalized to \RNP\ targets.
\medskip

Our proof is based on the notion of Markov convexity which was
introduced in \cite{LNP09}, with further important progress
achieved in \cite{MN13}.

\begin{definition}[\cite{LNP09}, we use a slightly modified version of \cite{MN13}]\label{D:MarkConv}
Let $\{X_t\}_{t\in \mathbb{Z}}$ be a Markov chain on a state space
$\Omega$. Given an integer $k\ge 0$,  we denote by $\{\widetilde
X_t(k)\}_{t\in \mathbb{Z}}$ the process which equals $X_t$ for
time $t\le k$, and evolves independently (with respect to the same
transition probabilities) for time $t > k$. Fix $p>0$. A metric
space $(X,d_X)$ is called {\it Markov $p$-convex with constant
$\Pi$} if for every Markov chain $\{X_t\}_{t\in \mathbb{Z}}$ on a
state space $\Omega$, and every $f : \Omega \to X$,
\begin{equation}\label{E:MarkConv}
\sum_{k=0}^{\infty}\sum_{t\in \mathbb{Z}}\frac{\mathbb{E}\left[
d_X\left(f(X_t),f\left(\widetilde
X_t\left(t-2^{k}\right)\right)\right)^p\right]}{2^{kp}} \le \Pi^p
\cdot \sum_{t\in \mathbb{Z}}\mathbb{E}
\big[d_X(f(X_t),f(X_{t-1}))^p\big].
\end{equation}
The least constant $\Pi$ for which~\eqref{E:MarkConv} holds for
all Markov chains is called the {\it Markov $p$-convexity
constant} of $X$, and is denoted $\Pi_p(X)$. We say that $(X,d_X)$
is {\it Markov $p$-convex} if $\Pi_p(X) < \infty$.
\end{definition}

Our proof of Theorem \ref{T:main} is based on the following
result:

\begin{theorem}\label{T:ThickNotMarkovC} A metric space with a thick family of
geodesics is not Markov $p$-convex for any $p\in(0,\infty)$.
\end{theorem}

Theorem \ref{T:ThickNotMarkovC} implies that thick families of
geodesics do not admit bilipschitz embeddings into the Heisenberg
group because it is known \cite[Theorem 7.4]{Li13} that the
Heisenberg group is Markov convex for some
$p\in(0,\infty)$.\medskip

Theorem \ref{T:ThickNotMarkovC} is a generalization of the result
of \cite[Section 3]{MN13} stating that the Laakso space (we mean
the Laakso space defined on \cite[p.~290]{LP01}) is not Markov
$p$-convex for any $p\in(0,\infty)$. It is easy to see that the
Laakso space has a thick family of geodesics.

\begin{remark} The Heisenberg group can be identified with $\mathbb{R}^3$ in such a way that all geodesics of the Heisenberg group with its subriemmanian metric
are spirals, projecting down to circles in two dimensions (see
\cite[Section 1.3]{Mon02}). With this representation it is easy to
verify that the family of $uv$-geodesics, where $u$ and $v$ are
elements of the Heisenberg group, is never thick. It is natural to
expect that one can prove Theorem 1.3 by combining this
description with some differentiability theory. I preferred to use
Markov convexity because I think that Theorem
\ref{T:ThickNotMarkovC} is of independent interest.
\end{remark}

\section{Proof of Theorem \ref{T:ThickNotMarkovC}}

Let $(M,d_M)$ be a metric space containing a thick family of
geodesics. We assume that each of the geodesics in the family has
length $1$ and is parameterized by the interval $[0,1]$.\medskip

The general idea of the proof is the same as the idea of the proof
of the fact that the Laakso space is not Markov convex in
\cite{MN13}. Namely, given $h\in \mathbb{N}$ we find in the thick
family of geodesics in $M$ (a thick family is necessarily
infinite) a finite collection $\mathbb{G}_h$ consisting of $2^h$
geodesics and a collection of finitely many points on each of
them, such that there is a Markov chain on this collection of
points with the the left-hand side of \eqref{E:MarkConv} greater
or equal than $C(p)h^{\frac1p}$ times the right-hand side of the
inequality \eqref{E:MarkConv} without $\Pi^p$.
\medskip

\noindent{\bf Short description of the Markov chain.} We introduce
the state space $\Omega$ as $\Omega=\mathbb{Z}\times
\mathbb{G}_h$. Let $\varphi$ be a positive integer (to be
specified later) and let $f:\Omega\to M$ be given by

\[
f(t,g)=\begin{cases} g(0) & \hbox{ if }t<0,\\
g(t2^{-\varphi}) &\hbox{ if }t\in \{0,1,2,3,\dots,2^\varphi\},\\
g(1) &\hbox{ if }t>2^\varphi.\end{cases}\]

The Markov chain $\{X_t\}_{t\in\mathbb{Z}}$ is defined as follows:

\begin{itemize}

\item $X_t=(t,g)$ for some $g\in\mathbb{G}_h$ (so the chain $X_t$
{\it remembers} the geodesic which {\it it is on}).

\item If $X_t=(t,g)$ and $t<0$ or $t\ge 2^\varphi$, then
$X_{t+1}=(t+1,g)$ with probability $1$.

\item If $X_t=(t,g)$ and $t\in \{0,1,2,3,\dots,2^\varphi-1\}$,
then $X_{t+1}=(t+1,\widehat g)$, where $\widehat g\in\mathbb{G}_h$
and either $\widehat g=g$ or $\widehat g=\widetilde g$, where
$\widetilde g$ is any geodesic of the family $\mathbb{G}_h$ which
has what we call a {\it crossing with $g$} in the interval
$[t2^{-\varphi}, (t+1)2^{-\varphi}]$. The probabilities of all
permissible choices of $\widehat g$ are the same. Crossings,
$\varphi$, and the family $\mathbb{G}_h$ of geodesics are defined
in such a way that a geodesic cannot have two crossings in one
interval of the form $[t2^{-\varphi}, (t+1)2^{-\varphi}]$.

\end{itemize}

We describe the needed notion of crossing below. At this point we
would like to mention that each crossing of geodesics corresponds
to their intersection, but not all of the intersections of
geodesics are crossings.
\medskip

The description of the allowed moves from one geodesic to another
in Theo\-rem~\ref{T:ThickNotMarkovC} is substantially more
complicated than in the case of the Laakso space in \cite{MN13},
because the geodesics can have infinitely many points of
intersection. Therefore to get the desired estimate we need the
Markov chain to move from one geodesic to another in a
well-organized manner, because we have lower estimates for
distances between geodesics only for small sets of pairs of points
(the only available estimate of this type is
\eqref{E:alphaSep}).\medskip

We label geodesics of $\mathbb{G}_h$ by {\it binary strings of
length $h$} and sets of crossings by vertices of a {\it binary
tree of depth $h-1$}.\medskip

Recall that a {\it binary tree $B_h$} of depth $h$ is a finite
graph whose vertices are finite sequences of $0$s and $1$s of
length at most $h$, including the empty sequence denoted
$\emptyset$; two vertices are joined by an edge if the
corresponding sequences are $(\theta_1,\dots,\theta_{n-1})$ and
$(\theta_1,\dots,\theta_{n-1},\theta_n)$ for some
$\theta_n\in\{0,1\}$ ($(\theta_1,\dots,\theta_{n-1})$ can be
empty).\medskip

We pick one element in the thick family of geodesics and label it
by the sequence consisting of $h$ zeros, so we denote it
$g_{(0,\dots,0)}$.  We apply the condition of Definition
\ref{D:ThickFam} to $g_{(0,\dots,0)}$ with control points $0$ and
$1$, and get a geodesic which we label $g_{(1,0,\dots,0)}$ and
points which we denote
$q_1^{\emptyset},\dots,q_m^{\emptyset}\in[0,1]$ and
$s^\emptyset_1,\dots s^\emptyset_{m+1}\in [0,1]$ such the
conditions of Definition \ref{D:ThickFam} are satisfied. We
introduce the set
\[R_\emptyset=\{0,q_1^{\emptyset},\dots,q_m^{\emptyset},1\}\]
This is the set of common crossings of all geodesics of
$\mathbb{G}_h$.\medskip

In the next step we pick two geodesics $g_{(0,1,0,\dots,0)}$ and
$g_{(1,1,0,\dots,0)}$ and find two subsets $R_{(0)}$ and $R_{(1)}$
of $[0,1]$. The sets $R_{(0)}$ and $R_{(1)}$ will be the sets of
common crossings of all geodesics whose labels start with $0$ and
with $1$, respectively.\medskip

To pick $g_{(0,1,0,\dots,0)}$ and $g_{(1,1,0,\dots,0)}$ we apply
the condition of Definition \ref{D:ThickFam} to $g_{(0,\dots,0)}$
and $g_{(1,0,\dots,0)}$, respectively, the collection of control
points defined as the union of two subsets:

\begin{itemize}

\item The points $q_1^{\emptyset},\dots,q_m^{\emptyset}$ and
$s^{\emptyset}_1,\dots, s^{\emptyset}_{m+1}\in [0,1]$.
\medskip

Observe that $g(0)$ is the same for all geodesics of the family,
this is the reason why we do not have to list $0$ among control
points. The same applies to $1$.\medskip

\item The points $k2^{-\gamma(1)}$, $k=1,\dots,2^{\gamma(1)}$,
where $\gamma(1)\in\mathbb{N}$ is sufficiently large. The
conditions on $\gamma(1)$ are the following:

\begin{enumerate}

\item $4m2^{-\gamma(1)}\le\frac{\alpha}{10}$, where $m$ is the
cardinality of $\{q_i^\emptyset\}$.

\item $2^{-\gamma(1)}\le \frac14
\min_i\left(s_i^{\emptyset}-q_{i-1}^{\emptyset}\right)$.

\end{enumerate}

\end{itemize}

We denote the sequences obtained by applying the condition of
Definition \ref{D:ThickFam} to $g_{(0,\dots,0)}$ by
$q_1^{(0)},\dots,q_{m(0)}^{(0)}\in[0,1]$ and $s^{(0)}_1,\dots,
s^{(0)}_{m(0)+1}\in [0,1]$, and the sequences obtained by applying
the condition of Definition \ref{D:ThickFam} to
$g_{(1,0,\dots,0)}$ by $q_1^{(1)},\dots,q_{m(1)}^{(1)}\in[0,1]$
and $s^{(1)}_1,\dots s^{(1)}_{m(1)+1}\in [0,1]$. We set
$R_{(0)}=\{q_1^{(0)},\dots,q_{m(0)}^{(0)}\}\backslash R_\emptyset$
and $R_{(1)}=\{q_1^{(1)},\dots,q_{m(1)}^{(1)}\}\backslash
R_\emptyset$. The set $R_{(0)}$ is the set of crossings of all
geodesics whose label starts with $0$. The set $R_{(1)}$ is the
set of crossings of all geodesics whose label starts with $1$.
\medskip

At this point we give a generic description which will be used for
all further selections of geodesics and sets of crossings.
\medskip

Suppose that we have already picked
$g_{(\theta_1,\dots,\theta_{n},0,\dots,0)}$ and constructed all
crossings sets $R_{(\theta_1,\dots,\theta_{k})}$, where
$\theta_1,\dots,\theta_{k}$ is an initial segment of
$\theta_1,\dots,\theta_{n-1}$, as well as the sequences
\begin{equation}\label{E:qGen} q_1^{(\theta_1,\dots,\theta_{n-1})},
\dots,q_{m(\theta_1,\dots,\theta_{n-1})}^{(\theta_1,\dots,\theta_{n-1})}\in[0,1]\end{equation}
and
\begin{equation}\label{E:sGen}s^{(\theta_1,\dots,\theta_{n-1})}_1,\dots,
s^{(\theta_1,\dots,\theta_{n-1})}_{m(\theta_1,\dots,\theta_{n-1})+1}\in
[0,1].\end{equation}

To pick the geodesic $g_{(\theta_1,\dots,\theta_{n},1,0,\dots,0)}$
we apply the conditions of Definition \ref{D:ThickFam} to
$g_{(\theta_1,\dots,\theta_{n},0,\dots,0)}$ and the collection of
control points defined as the union of two subsets:

\begin{itemize}

\item All points listed in \eqref{E:qGen} and \eqref{E:sGen}.

\item The points $k2^{-\gamma(n)}$, $k=1,\dots,2^{\gamma(n)}$
where $\gamma(n)$ is a sufficiently large number. The conditions
on $\gamma(n)$ are the following

\begin{equation}\label{E:GammaN}
4
m(\theta_1,\dots,\theta_{n-1})\cdot2^{-\gamma(n)}\le\frac{\alpha}{10}.
\end{equation}

\begin{equation}\label{E:GammaN2}
2^{-\gamma(n)}\le \frac14\min_{(\theta_1,\dots,\theta_{n-1})}
\min_i\left(s_i^{(\theta_1,\dots,\theta_{n-1})}-q_{i-1}^{(\theta_1,\dots,\theta_{n-1})}\right).
\end{equation}

\end{itemize}

We denote the obtained sequences

\begin{equation}\label{E:qGen+1} q_1^{(\theta_1,\dots,\theta_{n})},
\dots,q_{m(\theta_1,\dots,\theta_{n})}^{(\theta_1,\dots,\theta_{n})}\in[0,1]\end{equation}
and
\begin{equation}\label{E:sGen+1}s^{(\theta_1,\dots,\theta_{n})}_1,\dots,
s^{(\theta_1,\dots,\theta_{n})}_{m(\theta_1,\dots,\theta_{n})+1}\in
[0,1].\end{equation}

We introduce the set \[R_{(\theta_1,\dots,\theta_{n})}=
\left\{q_1^{(\theta_1,\dots,\theta_{n})},
\dots,q_{m(\theta_1,\dots,\theta_{n})}^{(\theta_1,\dots,\theta_{n})}\right\}\backslash\left(\bigcup_{k=0}^{n-1}R_{(\theta_1,\dots,\theta_{k})}\right),
\]
where by the set corresponding to $k=0$ we mean $R_\emptyset$. The
set $R_{(\theta_1,\dots,\theta_{n})}$ is the set of common
crossings of all geodesics whose label starts with
$(\theta_1,\dots,\theta_{n})$.\medskip

After we pick all geodesics for $\mathbb{G}_h$ and construct all
of the sets of crossings, we pick the number $\varphi\in
\mathbb{N}$. The choice of $\varphi$ should satisfy two
conditions:

\begin{itemize}

\item $2^{-\varphi}$ should be strictly less than the distance
between any two crossings which are crossings for the same
geodesic.

\item

\begin{equation}\label{E:ScndVarphi} 2^{-\varphi}\le\frac1{16}d_M\left(g_{(\theta_1,\dots,\theta_n,0,\dots,0)}(s_i^{(\theta_1,\dots,\theta_n)}),
g_{(\theta_1,\dots,\theta_n,1,0,\dots,0)}(s_i^{(\theta_1,\dots,\theta_n)})\right)\end{equation}
for all $(\theta_1,\dots,\theta_n)$ and $1\le i\le
m(\theta_1,\dots,\theta_n)+1$ satisfying
\begin{equation}\label{E:NE}
g_{(\theta_1,\dots,\theta_n,0,\dots,0)}(s_i^{(\theta_1,\dots,\theta_n)})\ne
g_{(\theta_1,\dots,\theta_n,1,0,\dots,0)}(s_i^{(\theta_1,\dots,\theta_n)})
\end{equation}

\end{itemize}

Now we are ready to complete the short description of the Markov
chain given at the beginning of the proof. Namely we provide more
details on the way in which Markov chain can move from one
geodesic to another. If $X_t=(t,
g_{(\theta_1,\dots,\theta_{h})})$, and the interval
$[t2^{-\varphi}, (t+1)2^{-\varphi}]$ contains a crossing labelled
by some initial segment $(\theta_1,\dots,\theta_d)$ of
$(\theta_1,\dots,\theta_h)$, then $X_{t+1}=(t+1,\widetilde g)$,
where $\widetilde g$ is any of the $2^{h-d}$ geodesics whose
labels have $(\theta_1,\dots,\theta_d)$ as their initial segment,
and each of these $2^{h-d}$ choices has the same probability.
Observe that the choice of $\varphi$ is such that a segment of the
form $[t2^{-\varphi}, (t+1)2^{-\varphi}]$ cannot contain more than
one crossing.
\medskip

For each collection $(\theta_1,\dots,\theta_n)$, $n<h$, we find a
subset $\left\{s_i^{(\theta_1,\dots,\theta_n)}\right\}_{i\in
A(\theta_1,\dots,\theta_n)}$ in the set
$\left\{s_i^{(\theta_1,\dots,\theta_n)}\right\}_{i=1}^{m(\theta_1,\dots,\theta_n)+1}$
which is sufficiently large in the sense that
\[\sum_{i\in
A(\theta_1,\dots,\theta_n)}d_M\left(g_{(\theta_1,\dots,\theta_n,0,\dots,0)}(s_i^{(\theta_1,\dots,\theta_n)}),
g_{(\theta_1,\dots,\theta_n,1,0,\dots,0)}(s_i^{(\theta_1,\dots,\theta_n)})\right)\ge\frac{\alpha}4.\]
We require that each $i\in A(\theta_1,\dots,\theta_n)$ satisfies
\eqref{E:NE} and two additional conditions needed for our
estimates; see conditions {\bf (a)} and {\bf (b)} below.\medskip

We estimate the sum in the left-hand side of \eqref{E:MarkConv}
from below as follows. We assign to each point
$s_i^{(\theta_1,\dots,\theta_n)}$ with $i\in
A{(\theta_1,\dots,\theta_n)}$ a scale $2^k$,
$(k=k(i,(\theta_1,\dots,\theta_n))\in\mathbb{N})$, an interval of
integers
\begin{equation}\label{E:Interval} I_{i,(\theta_1,\dots,\theta_n)}=[t,t+1,\dots,T],\end{equation} and a subset
$\mathbb{G}=\mathbb{G}_{(\theta_1,\dots,\theta_n)}$ in
$\mathbb{G}_h$, in such a way that no triple (scale, integer,
geodesic) is ever repeated. In the remainder of the argument we
use the following notation and terminology. If $\mathbb{G}$ is a
subset of $\mathbb{G}_h$ we write $X_{t}\in \mathbb{G}$ as a
shorthand for the condition $X_t=(t,g)$ with $g\in\mathbb{G}$. If
$X_t=(t,g)$, we say that $X_t$ is {\it on} $g$.

Then the left-hand side of \eqref{E:MarkConv} can be estimated
from below by

{\normalsize
\begin{equation}\label{E:DifSumOrd}\sum_{\stackrel{(\theta_1,\dots,\theta_n)}{i\in A(\theta_1,\dots,\theta_n)}}~\sum_{t\in
I_{i,(\theta_1,\dots,\theta_n)}}\frac{\left.\mathbb{E}\left[
d_X\left(f(X_t),f\left(\widetilde
X_t\left(t-2^{k}\right)\right)\right)^p\right|X_{t-2^k}\in
\mathbb{G}\right]\mathbb{P}\left(X_{t-2^k}\in
\mathbb{G}\right)}{2^{kp}},
\end{equation}}
where the conditional probability is with respect to the event
$X_{t-2^k}\in \mathbb{G}_{(\theta_1,\dots,\theta_n)}$. Note that
although it is not reflected in our notation, $k$ also depends on
$i,(\theta_1,\dots,\theta_n)$.\medskip

Now we describe how do we pick the scale $2^k$, the interval in
\eqref{E:Interval}, and the set of geodesics
$\mathbb{G}_{(\theta_1,\dots,\theta_n)}$.
\medskip

\noindent{\bf (1)} We pick $k$ to be the smallest positive integer
such that $2^k2^{-\varphi}$ exceeds
$s_i^{(\theta_1,\dots,\theta_n)}-q_{i-1}^{(\theta_1,\dots,\theta_n)}$
(we use $0$ instead of $q_{i-1}^{(\theta_1,\dots,\theta_n)}$ if
$i=1$).\medskip

\noindent{\bf (2)} We let
\begin{equation}\label{E:length}
L=L_{i,(\theta_1,\dots,\theta_n)}=d_M\left(g_{(\theta_1,\dots,\theta_n,0,\dots,0)}(s_i^{(\theta_1,\dots,\theta_n)}),
g_{(\theta_1,\dots,\theta_n,1,0,\dots,0)}(s_i^{(\theta_1,\dots,\theta_n)})\right)
\end{equation}
and introduce the interval \eqref{E:Interval} as the set of
$\tau\in\mathbb{Z}$ for which $\tau2^{-\varphi}$ is in the
interval of length $\frac14L$ which ends at the point
$s_i^{(\theta_1,\dots,\theta_n)}$.
\medskip

The set of such $\tau$ is nonempty because, by
\eqref{E:ScndVarphi}, $2^{-\varphi}<\frac1{16}L$. Furthermore,
\eqref{E:ScndVarphi} implies that
$2^{-\varphi}|I|\ge\frac18L$.\medskip

\noindent{\bf (3)} We define
$\mathbb{G}_{(\theta_1,\dots,\theta_n)}$ as the set of geodesics
whose labels start with $(\theta_1,\dots,\theta_n)$.
\medskip

Now we impose the second of the conditions under which
$i\in\{1,\dots,m(\theta_1,\dots,\theta_n)+1\}$ is included into
$A(\theta_1,\dots,\theta_n)$. (Below we introduce the third
condition which we label {\bf (b)}.)
\medskip

\noindent{\bf (a)} The interval of length
$\frac14L+2^k2^{-\varphi}$ which ends at the point
$s_i^{(\theta_1,\dots,\theta_n)}$ does not contain any crossings
belonging to $\bigcup_{k=0}^{n-1}R_{(\theta_1,\dots,\theta_k)}$.
\medskip

Observe that under the condition {\bf (a)} the conditional
expectation in \eqref{E:DifSumOrd} is at least
$\frac12\left(L-2\left(\frac14L\right)\right)^p=\frac1{2^{p+1}}L^p$.
The reason for this estimate is that the condition {\bf (a)}
implies that all of the geodesics in
$\mathbb{G}_{(\theta_1,\dots,\theta_n)}$ have
$q_{i-1}^{(\theta_1,\dots,\theta_n)}$ as their common crossing,
the crossing occurs ``after'' time $t-2^k$ if $t\in
I_{i,(\theta_1,\dots,\theta_n)}$, and there are no crossings which
could lead outside $\mathbb{G}_{(\theta_1,\dots,\theta_n)}$ in the
interval between $t-2^k$ and $t$ if $t\in
I_{i,(\theta_1,\dots,\theta_n)}$. Therefore with probability
$\frac12$ at the crossing $q_{i-1}^{(\theta_1,\dots,\theta_n)}$
the Markov chains $X_t$ and $\widetilde X_t$ will ``go'' in
different ``directions'',  one of them will ``go'' to
$g_{(\theta_1,\dots,\theta_n,0,\dots,0)}(s_i^{(\theta_1,\dots,\theta_n)})$,
and the other will ``go'' to
$g_{(\theta_1,\dots,\theta_n,1,0,\dots,0)}(s_i^{(\theta_1,\dots,\theta_n)})$.
It remains to use the triangle inequality.
\medskip

Next, it is easy to verify that the probability that $X_t$
$(t=0,1,\dots,2^\varphi)$  is on a geodesic $g$ is $2^{-h}$ if
$t2^{-\varphi}$ is not a crossing involving $g$. If
$t2^{-\varphi}$ is a crossing of $2^{h-n}$ geodesics, the
probability that $X_t$ is on one of them is $2^{-n}$. The
verification of this statement can be done by moving from $0$ to
$1$. Therefore the probability in \eqref{E:DifSumOrd} is $2^{-n}$.
\medskip

The third condition on  $i\in A(\theta_1,\dots,\theta_n)$
is\smallskip

\noindent{\bf (b)}
$L\ge\frac\alpha2\left(q_{i}^{(\theta_1,\dots,\theta_n)}-q_{i-1}^{(\theta_1,\dots,\theta_n)}\right)$,
where $L$ is defined in \eqref{E:length}.\medskip

Under the condition {\bf (b)} we can estimate each of the summands
in \eqref{E:DifSumOrd}. In fact, by the choice of $k$ we have
$2^k2^{-\varphi}<2\left(q_{i}^{(\theta_1,\dots,\theta_n)}-q_{i-1}^{(\theta_1,\dots,\theta_n)}\right)$.
Therefore $L>\frac\alpha4\cdot2^{k-\varphi}$ and
\[\frac{\frac1{2^{p+1}}L^p}{2^{kp}}>\frac{\alpha^p}{2^{3p+1}}2^{-p\varphi}.\]

Therefore, for $i$ satisfying the conditions {\bf (a)} and {\bf
(b)} each term in the sum \eqref{E:DifSumOrd} is $\ge
C2^{-n}2^{-p\varphi}$, where $C$ is a constant which depends only
on $\alpha$ and $p$.\medskip

Now we fix $(\theta_1,\dots,\theta_n)$ and consider the sum
\begin{small}
\begin{equation}\label{E:PartSum}\sum_{i\in A(\theta_1,\dots,\theta_n)}~\sum_{t\in
I_{i,(\theta_1,\dots,\theta_n)}}\frac{\left.\mathbb{E}\left[
d_X\left(f(X_t),f\left(\widetilde
X_t\left(t-2^{k}\right)\right)\right)^p\right|X_{t-2^k}\in
\mathbb{G}\right]\mathbb{P}\left(X_{t-2^k}\in
\mathbb{G}\right)}{2^{kp}}.
\end{equation}
\end{small}

As we observed above, the number of terms in the sum $\sum_{t\in
I_{i,(\theta_1,\dots,\theta_n)}}$ is at least
$2^\varphi\cdot\frac18L_{i,(\theta_1,\dots,\theta_n)}$. We shall
show that this implies that the sum in \eqref{E:PartSum} is at
least \begin{equation}\label{E:Des}\frac{2^\varphi}8\sum_{i\in
A(\theta_1,\dots,\theta_n)}L_{i,(\theta_1,\dots,\theta_n)}C2^{-n}2^{-p\varphi}\ge
\frac{C2^{(1-p)\varphi}2^{-n}}8\left(\frac\alpha2-\frac\alpha{10}\right).\end{equation}

To get this we used the inequality \[\sum_{i{\rm
~satisfies~}{\bf(b)}}L_{i,(\theta_1,\dots,\theta_n)}\ge\frac\alpha2,\]
which follows from \eqref{E:alphaSep} by the Markov-type
inequality. The condition \eqref{E:NE} is included in {\bf (b)},
but we also have to exclude $i$ which fail to satisfy {\bf (a)}.
\medskip

Comparing condition {\bf (a)} with our definitions we see that it
suffices to require that $s_i^{(\theta_1,\dots,\theta_n)}$ is not
in the interval of length
$2(s_i^{(\theta_1,\dots,\theta_n)}-q_{i-1}^{(\theta_1,\dots,\theta_n)})+\frac14L\le
3(s_i^{(\theta_1,\dots,\theta_n)}-q_{i-1}^{(\theta_1,\dots,\theta_n)}\le
3\cdot 2^{-\gamma(n)}$ (see \eqref{E:GammaN}) following one of the
elements of $\bigcup_{k=0}^{n-1}R_{(\theta_1,\dots,\theta_k)}$.
Therefore the  total length of the intervals
$[q_{i-1}^{(\theta_1,\dots,\theta_n)},q_{i}^{(\theta_1,\dots,\theta_n)}]$
which have to be excluded does not exceed
$4m(\theta_1,\dots,\theta_{n-1})\cdot
2^{-\gamma(n)}\le\frac{\alpha}{10}$ (see \eqref{E:GammaN}). It is
clear that the sum of $L_{i,(\theta_1,\dots,\theta_{n})}$ over all
excluded in this way $i$ also does not exceed $\frac{\alpha}{10}$.
The inequality \eqref{E:Des} follows.\medskip

Therefore the sum in \eqref{E:PartSum} is $\ge
2^{(1-p)\varphi}2^{-n}C(\alpha,p)$. Adding over
$(\theta_1,\dots,\theta_n)$ for fixed $n$ we get $\ge
2^{(1-p)\varphi}C(\alpha,p)$. Adding over $n=0,1,\dots,h-1$ we get
$\ge 2^{(1-p)\varphi}C(\alpha,p)h$. On the other hand, the sum in
the right-hand side of \eqref{E:MarkConv} is $2^{-p\varphi}\cdot
2^\varphi=2^{(1-p)\varphi}$. We get $\Pi^p\ge C(\alpha,p)h$, which
is the desired inequality.\medskip

To complete the proof we need to explain why for different choices
of $i$ and $(\theta_1,\dots,\theta_n)$ the sets of triples (scale,
integer, geodesic) are disjoint. First let us consider the case
where $(\theta_1,\dots,\theta_{n_1})$ and
$(\tilde\theta_1,\dots,\tilde\theta_{n_2})$ are such that $n_1\ne
n_2$. Observe that $2^{k(i,(\theta_1,\dots,\theta_{n_1}))}$ is
$2$-equivalent to
\begin{equation}\label{E:i}2^\varphi(s_i^{(\theta_1,\dots,\theta_{n_1})}-q_{i-1}^{(\theta_1,\dots,\theta_{n_1})}),\end{equation}
and $2^{k(j,(\tilde\theta_1,\dots,\tilde\theta_{n_2}))}$ is
$2$-equivalent to
\begin{equation}\label{E:j}2^\varphi(s_j^{(\tilde\theta_1,\dots,\tilde\theta_{n_2})}-q_{j-1}^{(\tilde\theta_1,\dots,\tilde\theta_{n_2})}),\end{equation}
and for $n_1\ne n_2$ the numbers \eqref{E:i} and \eqref{E:j}
cannot be $4$-equivalent, see \eqref{E:GammaN2} and the
description of the choice of geodesics
$g_{(\theta_1,\dots,\theta_h)}$.\medskip

If $n_1=n_2$, but $(\theta_1,\dots,\theta_{n_1})$ is not the same
as $(\tilde\theta_1,\dots,\tilde\theta_{n_2})$, then the families
$\mathbb{G}_{i,(\theta_1,\dots,\theta_{n_1})}$ and $\mathbb{G}_{j,
(\tilde\theta_1,\dots,\tilde\theta_{n_2})}$ do not contain common
geodesics.\medskip

Finally, if we consider labels $i,(\theta_1,\dots,\theta_n)$ and
$j,(\theta_1,\dots,\theta_n)$, then either
\[k(i,(\theta_1,\dots,\theta_n))\ne
k(j,(\theta_1,\dots,\theta_n))\] (and we are done) or
\[k(i,(\theta_1,\dots,\theta_n))=
k(j,(\theta_1,\dots,\theta_n)).\] In the latter case, as is easy
to check, the intervals $I_{i,(\theta_1,\dots,\theta_n)}$ and
$I_{j,(\theta_1,\dots,\theta_n)}$ (see {\bf (2)} for the
definition) are disjoint.

\end{large}

\begin{small}

\renewcommand{\refname}{
\end{small}

\end{document}